# Huellas del caminar de los números: una tarea de visualización dinámica para la comprensión de los números decimales en educación secundaria

**Felix De la Cruz Serrano**
I. E. María Reiche
Perú
feldese@gmail.com

**Resumen**

El estudio de los números decimales en la educación secundaria suele abordarse desde enfoques algorítmicos, lo que limita la comprensión de su estructura. En este trabajo se presenta la tarea Huellas del caminar de los números, una propuesta de visualización dinámica orientada a favorecer la comprensión de los números decimales a partir de la exploración de sus desarrollos infinitos.

La tarea se basa en la asignación de vectores a las cifras decimales del 0 al 9, de modo que la secuencia de cifras genera un recorrido geométrico dinámico en el plano. Mediante el uso de GeoGebra como entorno de visualización, los estudiantes pueden observar, comparar e interpretar huellas asociadas a distintos tipos de números, como decimales exactos, periódicos e irracionales, identificando regularidades visuales vinculadas a su comportamiento decimal.

El análisis se desarrolla desde una perspectiva teórica–didáctica, utilizando el Espacio de Trabajo Matemático como lente interpretativa para caracterizar el potencial del diseño de la tarea. Asimismo, se justifica el uso puntual de herramientas de IA generativa exclusivamente como apoyo instrumental al cálculo, sin desplazar el foco del razonamiento matemático.

*Palabras clave*: visualización matemática, números decimales, GeoGebra, diseño de tareas.

**Abstract**

The study of decimal numbers in secondary education is often approached from algorithmic perspectives, which limits the understanding of their structure. This paper presents the task Huellas del caminar de los números, a dynamic visualization proposal aimed at supporting the understanding of decimal numbers through the exploration of their infinite decimal expansions.

The task is based on assigning vectors to the decimal digits from 0 to 9, so that the sequence of digits generates a dynamic geometric path in the plane. Through the use of GeoGebra as a visualization environment, students can observe, compare, and interpret traces associated with different types of numbers, such as terminating decimals, repeating decimals, and irrational numbers, identifying visual regularities linked to their decimal behavior.

The analysis is developed from a theoretical-didactical perspective, using the Mathematical Working Space as an interpretative lens to characterize the potential of the task design. Likewise, the paper justifies the punctual use of generative AI tools exclusively as instrumental support for computation, without shifting the focus away from mathematical reasoning.

*Keywords:* mathematical visualization, decimal numbers, GeoGebra, task design.

## Introducción

En la enseñanza secundaria, los números decimales suelen presentarse como extensiones operativas de las fracciones, enfatizando procedimientos de cálculo, transformaciones simbólicas y reglas de clasificación. Este enfoque, si bien necesario, con frecuencia deja en segundo plano la comprensión de la estructura del desarrollo decimal y las propiedades que permiten distinguir distintos tipos de números reales. Como consecuencia, los estudiantes tienden a percibir los decimales como secuencias de cifras sin significado matemático profundo.

Diversos estudios han señalado el papel central de la visualización en la construcción de significado matemático, particularmente cuando se trata de objetos abstractos cuya comprensión requiere coordinar múltiples representaciones (Duval, 2006). Desde esta perspectiva, la visualización no se concibe como un recurso ilustrativo accesorio, sino como un proceso cognitivo que permite hacer visibles propiedades y estructuras que permanecen ocultas en tratamientos exclusivamente simbólicos (Arcavi, 1999). En este sentido, los entornos de geometría dinámica ofrecen oportunidades para explorar relaciones, regularidades y patrones que difícilmente emergen en enfoques centrados únicamente en el cálculo.

En este trabajo se presenta la tarea Huellas del caminar de los números, diseñada para promover una aproximación visual a los números decimales mediante el uso de GeoGebra. El propósito de la comunicación es analizar el potencial didáctico del diseño de esta tarea desde una perspectiva teórica–didáctica, poniendo el foco en la visualización como eje del trabajo matemático y en el rol del entorno dinámico como mediador de la actividad del estudiante.

## Marco teórico

### Visualización y Matemáticas Visibles

La visualización en educación matemática no se limita a la producción de imágenes ilustrativas, sino que constituye un medio para pensar, explorar y construir significado matemático. Visualizar implica coordinar acciones, representaciones y razonamientos para acceder a propiedades que no son evidentes de manera inmediata en el registro simbólico (Duval, 2006). En esta línea, Arcavi (1999) caracteriza la visualización como una capacidad, un proceso y un producto del pensamiento matemático, subrayando su papel en la posibilidad de "ver lo invisible", es decir, de hacer accesibles estructuras matemáticas que no se manifiestan directamente en las expresiones formales.

Desde este enfoque, la visualización adquiere un rol central en la actividad matemática escolar, favoreciendo el diseño de tareas en las que los estudiantes puedan observar, manipular y analizar representaciones que hagan visibles aspectos estructurales de los objetos matemáticos. En el caso de los números decimales, este tipo de tareas permite desplazar la atención desde el cálculo operativo hacia la identificación de regularidades, patrones y comportamientos asociados a su desarrollo infinito, promoviendo una comprensión más profunda de su naturaleza matemática.

**El Espacio de Trabajo Matemático como lente interpretativa**

Para analizar el potencial del diseño de la tarea, se adopta el Espacio de Trabajo Matemático (ETM) como marco interpretativo. Desde esta perspectiva, el trabajo matemático escolar se concibe como la articulación entre acciones instrumentales, interpretaciones semióticas y producciones discursivas, mediadas por el uso de artefactos (Kuzniak, 2011). Este enfoque permite dar cuenta de cómo la actividad matemática se organiza a partir de la coordinación entre herramientas, representaciones y razonamientos.

En coherencia con investigaciones sobre el uso de tecnologías digitales en educación matemática, el ETM ofrece una estructura adecuada para analizar propuestas en las que la tecnología se integra como un medio al servicio del trabajo matemático, y no como un fin en sí mismo. En este sentido, Drijvers (2013) señala que el potencial didáctico de las herramientas digitales depende fundamentalmente del diseño de las tareas, del rol del docente y del contexto educativo en el que se insertan. Desde esta perspectiva, el uso de GeoGebra en la tarea de visualización propuesta se concibe como un artefacto mediador que favorece la coordinación entre diferentes registros de representación y la construcción de significado matemático.

En este estudio, el ETM no se emplea con fines de análisis empírico exhaustivo, sino como una lente conceptual que permite caracterizar cómo el diseño de la tarea de visualización promueve dicha coordinación y organiza el trabajo matemático en estudiantes de educación secundaria.

## Metodología y diseño de la tarea Huellas del caminar de los números

**Intencionalidad didáctica**

La tarea se diseña con el propósito de favorecer una comprensión estructural de los números decimales a partir de la visualización de sus desarrollos infinitos. En lugar de centrar la actividad en el cálculo, se propone que los estudiantes exploren representaciones visuales dinámicas que les permitan reconocer regularidades asociadas a la periodicidad decimal y distinguir comportamientos propios de números exactos, periódicos e irracionales.

La visualización se concibe como un medio para formular conjeturas, comparar casos y argumentar a partir de lo observado, promoviendo un trabajo matemático que articule lo visual con lo conceptual.

**Descripción matemática de la tarea**

La tarea se basa en la asignación de vectores a cada una de las cifras decimales del 0 al 9. Dado un número real expresado en su desarrollo decimal, cada cifra se asocia a un vector fijo en el plano. La secuencia de cifras genera así una concatenación de desplazamientos cuya representación gráfica constituye una trayectoria o "huella".

Formalmente, si a cada cifra $d \in \{0; 1; 2; 3; 4; 5; 6; 7; 8; 9\}$ se le asigna un vector $\vec{v}_d$, el recorrido se obtiene sumando sucesivamente los vectores correspondientes a las cifras del desarrollo decimal del número. El rastro resultante puede observarse de manera estática o animada, permitiendo analizar su forma global y sus regularidades locales. La Figura 1 ilustra una huella generada a partir del desarrollo decimal del número racional 1/14, en la que es posible identificar un patrón visual regular asociado a la periodicidad del desarrollo decimal.

**Figura 1**
*Huella generada a partir del desarrollo decimal del número racional 1/14 mediante la asignación vectorial de cifras en GeoGebra.*

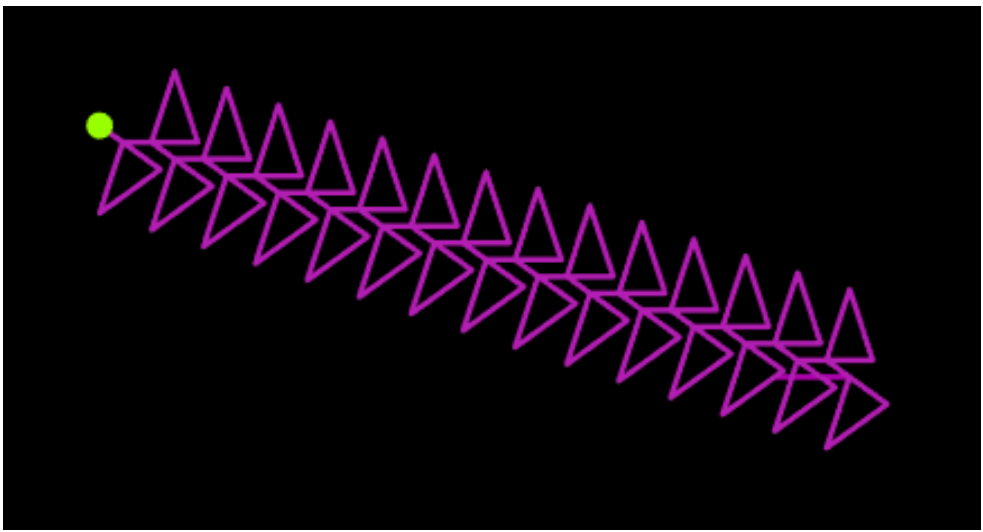

Este procedimiento permite obtener representaciones visuales diferenciadas según la naturaleza del número considerado. En los decimales exactos, el recorrido se detiene tras un número finito de pasos; en los decimales periódicos, emergen patrones repetitivos; mientras que en los números irracionales no se observan regularidades evidentes de repetición.

**Uso de GeoGebra y herramientas digitales con IA generativa**

La tarea se implementa en GeoGebra, que actúa como entorno de visualización dinámica y permite construir, animar y comparar las huellas generadas. Las funcionalidades del software facilitan la observación progresiva del recorrido y la exploración de distintos casos, situando la visualización como eje del trabajo matemático.

Una limitación técnica de GeoGebra es la restricción en el número de cifras decimales que se generan directamente al dividir fracciones, lo cual resulta insuficiente cuando se busca analizar patrones visuales asociados a desarrollos decimales extensos. Para superar esta limitación, se recurre de manera puntual a herramientas de IA generativa con el fin de obtener listas ampliadas de cifras decimales en un formato compatible con GeoGebra.

En particular, se diseñó un procedimiento automatizado basado en un modelo de lenguaje que, a partir de un prompt específico, genera de manera controlada las cifras decimales correspondientes a un número dado y las entrega como una lista utilizable en el entorno de visualización. La IA se emplea exclusivamente como apoyo instrumental al cálculo, sin intervenir en la interpretación matemática ni en la toma de decisiones por parte del estudiante, manteniéndose el foco en el análisis visual y el razonamiento matemático.

## Resultados: análisis del potencial didáctico desde el ETM

El análisis del potencial didáctico de la tarea se realiza desde el ETM, considerando las posibilidades de trabajo matemático que el diseño habilitó.

### Acciones instrumentales

GeoGebra permite acciones como la animación del recorrido, la observación paso a paso del rastro y la comparación entre huellas asociadas a distintos números. Estas acciones favorecen procesos de anticipación, control visual y verificación de conjeturas, organizando la actividad matemática en torno a la exploración.

### Articulación de registros semióticos

La tarea exige coordinar el registro numérico simbólico del desarrollo decimal con un registro geométrico dinámico. Esta conversión semiótica obliga a interpretar visualmente la estructura del número y a relacionar regularidades gráficas con propiedades del desarrollo decimal. Mientras que en el caso de números racionales periódicos es posible identificar patrones visuales asociados a la repetición del período, la huella correspondiente al número irracional $\pi$ no presenta regularidades de repetición, generando un recorrido más irregular y no estructurado, coherente con la ausencia de periodicidad en su desarrollo decimal, como se muestra en la Figura 2.

### Producción de discurso matemático

A partir de la visualización, los estudiantes pueden formular afirmaciones y justificar visualmente la presencia o ausencia de periodicidad, transitando desde explicaciones basadas en la observación hacia nociones matemáticas progresivamente más formales.

**Figura 2**

*Huella generada a partir del desarrollo decimal del número irracional π mediante la asignación vectorial de cifras en GeoGebra.*

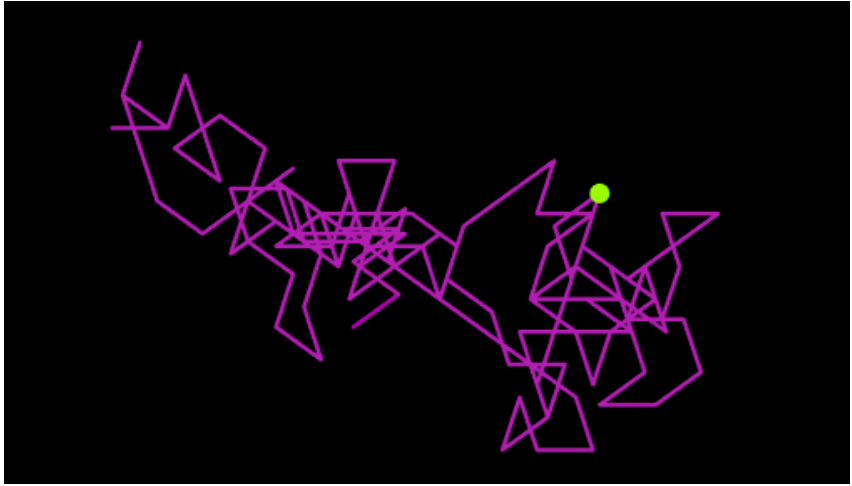

### Discusión

El análisis del diseño de la tarea Huellas del caminar de los números pone de manifiesto el papel de la visualización dinámica como eje del trabajo matemático en educación secundaria. La articulación entre acciones instrumentales en GeoGebra, la conversión entre registros numéricos y geométricos y la producción de discurso matemático permite desplazar el foco desde procedimientos algorítmicos hacia la exploración de regularidades y patrones asociados al desarrollo decimal.

Desde la perspectiva del Espacio de Trabajo Matemático, la tarea favorece una coordinación coherente entre el uso del artefacto digital, las interpretaciones semióticas y las explicaciones matemáticas basadas en la observación. En este sentido, la visualización no actúa como un recurso ilustrativo, sino como un medio para la construcción de significado, alineado con el enfoque de Matemáticas Visibles.

Asimismo, el uso de GeoGebra como entorno de visualización dinámica resulta clave para estructurar la actividad matemática, mientras que la incorporación puntual de herramientas de IA generativa se mantiene en un plano estrictamente instrumental, ampliando las posibilidades de exploración sin desplazar el protagonismo del razonamiento matemático del estudiante.

### Conclusión

En esta comunicación se ha presentado y analizado la tarea Huellas del caminar de los números como una propuesta de visualización dinámica orientada a favorecer la comprensión de los números decimales en educación secundaria. El análisis realizado muestra que el diseño de la tarea permite abordar propiedades estructurales del desarrollo decimal a partir de experiencias visuales que promueven la exploración, la formulación de conjeturas y la argumentación matemática.

El uso de GeoGebra como entorno de visualización se consolida como un elemento central del diseño, al posibilitar la construcción y comparación de huellas asociadas a distintos tipos de números. Por su parte, la utilización puntual de herramientas de IA generativa se justifica como apoyo instrumental al cálculo, manteniendo el foco en el trabajo matemático y evitando que la tecnología adquiera un protagonismo que desplace el razonamiento matemático.

Finalmente, este trabajo abre la posibilidad de futuras investigaciones centradas en el análisis de producciones estudiantiles durante la implementación de la tarea en el aula, así como en la extensión del enfoque visual a otros objetos matemáticos relevantes en la educación secundaria.

## Referencias y bibliografía